\documentclass[preprint,12pt,authoryear]{elsarticle}

\usepackage{amssymb}
\usepackage{amsthm, amsmath}
\usepackage{dsfont}
\usepackage{graphicx}
\usepackage{color}

\journal{Special Issue SPA}

\newcommand\R{\mathbb{R}}
\newcommand\E{\mathbb{E}}
\renewcommand\P{\mathbb{P}}
\newcommand\1{\mathds{1}}
\renewcommand\H{{\mathcal{H}}}
\newcommand\var{\mathrm{Var } }
\newcommand\Crit{\mathrm{Crit } }

\newtheorem{prop}{Proposition}
\newtheorem{lemma}[prop]{Lemme}
\newtheorem{thm}[prop]{Theorem}

\begin{document}

\begin{frontmatter}






\title{Minimal penalty for Goldenshluger-Lepski method}


\author[laborsay]{Lacour, C.}
\author[laborsay]{Massart, P.}
\address[laborsay]{Laboratoire de Math\'ematiques d'Orsay, Universit\'e Paris-Sud}

\begin{abstract}
This paper is concerned with adaptive nonparametric estimation using the Goldenshluger-Lepski selection method. This estimator selection method is based on pairwise comparisons between estimators with respect to some loss function. The method also involves a penalty term that typically needs to be  large enough in order that the method works (in the sense that one can prove some oracle type inequality for the selected estimator). In the case of density estimation with kernel estimators and a quadratic loss, we show that the procedure fails if the penalty term is chosen smaller than some critical value for the penalty: the minimal penalty. More precisely we show that the quadratic risk of the selected estimator explodes when the penalty is below this critical value while it stays under control when the penalty is above this critical value. This kind of phase transition phenomenon for penalty calibration has already been observed and proved for penalized model selection methods in various contexts but appears here for the first time for the Goldenshluger-Lepski pairwise comparison method.
Some simulations illustrate the theoretical results and lead to some hints on how to use the theory to calibrate the method in practice.
\end{abstract}


\end{frontmatter}

\section{Introduction}

Adaptive estimation is a challenging task in nonparametric estimation. Many methods have been proposed and studied in the literature. Most of them rely on some data-driven selection of an estimator among a given collection. Wavelet thresholding \citep{DJKP}, Lepski's method \citep{Lepski90}, and model selection \citep*{BBM} (see also \cite{Birge01} for the link between model  selection 
and Lepski's method) belong to this category. Designing proper estimator selection is an issue by itself.
From a constructive point of view, it is a crucial step towards adaptive estimation. For instance, selecting a bandwidth for kernel estimators in density estimation means that you are able to estimate the density without specifying its degree of smoothness in advance. 
Recently an interesting new estimator selection procedure has been introduced by  \cite{GL08}. 
Assume that one wants to estimate some unknown function $f$ belonging to some function space endowed with some norm $\|.\|$. Assume also that we have at our disposal some collection of estimators $(\hat f_h)_{h\in \H}$ indexed by some parameter $h$, the issue being to select some estimator $\hat f_{\hat h}$ among this collection.
The Goldenshluger-Lepski method proposes to select $\hat h$ as a minimizer 
of $B(h)+V(h)$ with $$B(h)=\sup\{[\|\hat f_{h'}- \hat f_{h,h'}\|^2-V(h')]_+, \: h'\in \H\}$$
where $x_+$ denotes the positive part $\max(x,0)$ and where $\hat f_{h,h'}$ are auxiliary (typically oversmoothed) estimators and 
$V(h)$ is a penalty term (called "majorant" by Goldenshluger and Lepski) to be suitably chosen.
They first develop their methodology in the white noise framework \citep{GL08,GL09}, next for density estimation \citep{GL10} and then for various other frameworks \citep{GL13}. 
Their initial motivation was to provide adaptive procedures for multivariate and anisotropic estimation and they used the versatility of their method to prove that the selected estimators can achieve minimax rates of convergence over some very general classes of smooth functions  \cite[see][]{GL14}. 
To this purpose, they have established oracle inequalities to ensure that, if $V(h)$ is well chosen, 
the final estimator $\hat f_{\hat h}$ is almost as efficient as the best one in the collection.
The Goldenshluger-Lepski methodology has already been fruitfully applied in various contexts: 
transport-fragmentation equations \citep{doumic2012}, anisotropic deconvolution \citep{comtelacour13}, warped bases regression \citep{chagny2013} among others
(see also \cite{BLR} which contains some explanation on the methodology).
We cannot close this paragraph without citing the nice work of \cite{LLP08}, who have independently introduced a very similar method, in order
to adapt the model selection point of view to pointwise estimation.

In this paper we focus on the issue of calibrating the penalty term $V$. As we mentioned above the "positive" known results are of the following kind:  the method performs well (at least from a theoretical view point) when $V$ is well chosen. More precisely one is able to prove oracle inequalities only if $V$ is not too small. But the issue is now: what is the minimal (or the optimal) value for $V$ to preserve (or optimize) the performance of the method? Here we consider this issue from a theoretical point of view but actually it is a crucial issue for a practical implementation of the method. In this paper we focus on the (simple) classical bandwidth selection issue for kernel estimators in the framework of univariate density estimation.
The main contribution of this paper is to highlight a phase transition phenomenon that can be roughly described as follows. For some critical quantity $V_0$ (that we call "minimal penalty") if the penalty term $V$ is defined as $V=aV_0$ then either $a<1$ and the risk 
$\E\|f-\hat f_{\hat h}\|^2$ is proven to be dramatically suboptimal, or $a>1$ and the risk remains under control. This kind of phase transition phenomenon and its possible use for penalty calibration appeared for the first time in \cite{birgemassart07} in the context of Gaussian penalized model selection.  It is interesting to see that the same phenomenon occurs for a pairwise comparison based selection method such as the Goldenshluger-Lepski method.

Proofs are extensively based on concentration inequalities. In particular, left tail concentration inequalities are used to prove the explosion result below the critical value for the penalty. Although the probabilistic tools are non asymptotic by essence, they merely allow us to justify that suprema of empirical processes are well concentrated around  their expectations and the approximations that we make on those expectations are indeed asymptotic. Needless to say this means that our final results are (unfortunately) a bit of an asymptotic nature, at least as far as the identification of the critical value $a=1$ is concerned. To be more concrete, we mean that for a given unknown density and a given sample size $n$, it is unclear that a phase transition phenomenon (if any) should occur at the critical value  $a=1$ as predicted by the (asymptotic) theory.  But still, because of the concentration phenomenon, one can hope that some phase transition does occur (even non asymptotically) at some critical value even though it is not equal (or even close) to the (asymptotic) value $a=1$.
To check this, we have also implemented numerical simulations. These simulations allow us to understand what should be retained from the theory as a typical behavior of the method. In fact the simulations confirm the above scenario. It turns out that the phase transition does occur when you run simulations even though the critical point is not located at $a=1$. This is actually what should be retained from the theory (at least from our point of view). The fact that some phase transition does occur is good news for the calibration issue because this means that in practice you can detect the critical value from the data (forgetting about the asymptotic value $a=1$). Then you can hope to use this value to elaborate some fully data-driven and non asymptotic calibration of the method.  We conclude the paper with providing some hints on how to perform that explicitly.

In Section~\ref{sec:sup} we specify the statistical framework and we recall the oracle inequality that can be obtained in the framework of density estimation. Then Section~\ref{minipen} contains our main theorem about minimal penalty. This result is illustrated by some simulations (Section~\ref{simus}).
Finally, some proofs are gathered in Section~\ref{preuves} after some concluding remarks.

\section{Kernel density estimation framework and upper bound on the risk}\label{sec:sup}

We consider independent and identically distributed real variables $X_1, \dots, X_n$ with unknown density $f$ with respect to the Lebesgue measure on the real line. Let $\|.\|$ denote the $L^2$ norm with respect to the Lebesgue measure. 
For each positive number $h$ (the bandwidth) we can define the classical kernel density estimator 
\begin{equation}\label{classicest}
\hat f_{h}(x)=\frac 1n \sum_{i=1}^n K_h(x-X_i)
\end{equation}
where $K$ is a kernel and $K_h=K(./h)/h$. 
We assume here that the function to be estimated is univariate and we study the Goldenshluger-Lepski methodology without oversmoothing. This means that 
we do not use auxiliary estimators. We could actually prove the same results for the original method but the proofs are more involved and we decided to keep the proofs as simple as possible trying not to hide the heart of the matter.

To be more precise the procedure that we study is the following one: 
%
%
starting from some (finite) collection of estimators $\{\hat f_{h}, h\in \H\}$, 
we set
\begin{equation}\label{defBV}
B(h)=\sup_{h'\leq h}\left[\|\hat f_{h'}-\hat f_{h}\|^2-V(h')\right]_+
\qquad \text{ with }V(h')=a\frac{\|K_{h'}\|^2}n
\end{equation}
with  $a$ being the tuning parameter of interest.
Then the selected bandwidth is defined by
\begin{equation}\label{defhath}
\hat{h}=\arg\min_{h\in\mathcal{H}}\left\{ B(h)+V(h)\right\}. 
\end{equation}
It is worth noticing that the penalty term $V(h)$ which is used here is exactly proportional to the integrated variance of the corresponding estimator.

We introduce the following notation:
\begin{eqnarray*}
&& f_h:=\E(\hat f_h), \quad  h_{\min}:=\min{\H}, \quad h_{\max}:=\max{\H} \\
&& D(h):=\max(\sup_{h'\leq h}\|f_{h'}-f_h\|,\|f-f_h\|)\leq 2\sup_{h'\leq h}\|f_{h'}-f\|
\end{eqnarray*}

We assume that the kernel verifies assumption
\begin{description}
\item[(K0)]
$\int |K|=1$, $\|K\|<\infty$  and 
\begin{equation*}
\forall \;0\leq x\leq 1 \qquad \frac{\langle K, K(x.) \rangle}{\|K\|^2}\geq 1.
\end{equation*}
\end{description}
Assumption {\bf (K0)} is satisfied whenever the kernel $K$ is nonnegative and unimodal with a mode at 0. Indeed in this case $K(xu) \geq K(u)$ for all $u\in \R$ and $x\in[0,1]$.
This is verified for classical kernels (Gaussian kernel, rectangular kernel, Epanechnikov kernel, biweight kernel; see Lemma~\ref{noyaux}).
This entails that for all $h'\leq h$, $\|K_{h'}-K_h\|^2\leq \|K_{h'}\|^2-\|K_h\|^2$. This Pythagore type inequality is a one of the key properties that we shall use for proving our results. 

Let us now recall the positive results that can be obtained for the selection method if $a$ is well chosen. 

\begin{prop} \label{oi} Assume that  $f$ is bounded and $K$ verifies {\bf (K0)}. Let $\hat f_{\hat h}$ be the selected estimator defined by \eqref{classicest}, \eqref{defBV}, \eqref{defhath}. Assume that the parameter $a$ in the penalty $V$ satisfies $a>1$. 

$\bullet$
There exist some positive constants $C_0>0$ and $c>0$ such that, with probability larger than
 $$1-2\sum_{h\in \H}\sum_{h'\leq h}\max(e^{-c\sqrt {n}},e^{-c/h'}),$$ 
the following holds
$$\|\hat f_{\hat h}-f\|\leq  C_0\inf_{h\in\mathcal H}\left\{   D(h)+\sqrt{a}\frac{\|K_h\|}{\sqrt{n}}\right\}.$$
The values $C_0=1+\sqrt{2(1+({a}^{1/3}-1)^{-1})}$ and $c=({a}^{1/3}-1)^2\min(\frac1{24},
\frac{\|K\|^2}{6\|f\|_\infty})$ are suitable.

$\bullet$
Moreover, if $n^{-1}\leq h_{\min}\leq h_{\max}\leq\log^{-2}(n)$, 
there exists a positive constant $C$ depending only on $\|K\|$ and $\|f\|_\infty$ such that 
$$\E\|\hat f_{\hat h}-f\|^2\leq 2C_0^2\inf_{h\in\mathcal H}\left\{   D^2(h)+a\frac{\|K_h\|^2}{{n}}\right\}+Cn{|\H|^2}e^{-\frac{({a}^{1/3}-1)^2}{C}\log^2(n)}$$
($C=\max(24,6\|f\|_\infty/\|K\|^2,4\|f\|_\infty+4\|K\|^2)$ works).
\end{prop}

We recognize in the right-hand side of the oracle type inequalities above the classical bias variance tradeoff. This oracle inequality shows that the Goldenshluger-Lepski methodology works when $a>1$,
at least for $n$ larger than some integer depending on $a$ and the true density. 
From a non asymptotic perspective this "positive result" should be understood with caution, it is clear from the analysis of the behavior of the constants involved with respect to $a$ that these constants are worse when $a$ is close to 1.

The proof of Proposition~\ref{oi} is postponed in  Section~\ref{preuveoi}. It is based on the following concentration result (adapted from \cite{KleinRio}) and more precisely on inequality~\eqref{concsup} below. 
\begin{lemma}\label{conc}
Let $X_1,\dots,X_n$ be a sequence of i.i.d. variables and 
$\nu(t)=n^{-1}\sum_{i=1}^n [g_t(X_i)-{\mathbb E}(g_t(X_i))]$
for $t$ belonging to a countable set of functions ${\mathcal F}$.
Assume that for all $t\in {\mathcal F}$ $\|g_t\|_\infty\leq b$ and
$ \var (g_t(X_1))\leq v$.
Denote  $H=\E(\sup_{t\in {\mathcal F}}\nu(t))$.
Then, for any $\varepsilon>0$, for $H'\geq H$,
\begin{eqnarray}
 \P(\sup_{t\in {\mathcal F}}\nu(t)\geq (1+\varepsilon)H')\leq \max\left(\exp\left(-\frac{\varepsilon^2}{6}\frac{nH'^2}{v}\right),
\exp\left(-\frac{\min(\varepsilon,1)\varepsilon}{24}\frac{n H'}{b}\right)\right) \label{concsup}\\
\P(\sup_{t\in {\mathcal F}}\nu(t)\leq H-\varepsilon H')\leq \max\left(\exp\left(-\frac{\varepsilon^2}{6}\frac{nH'^2}{v}\right),
\exp\left(-\frac{\min(\varepsilon,1)\varepsilon}{24}\frac{n H'}{b}\right)\right)\label{concinf}
\end{eqnarray}
Moreover
\begin{equation}\label{varsup}
\var (\sup_{t\in {\mathcal F}}\nu(t))\leq \frac vn+4\frac{bH}n
\end{equation}
\end{lemma}

\section{Minimal penalty}\label{minipen}

In this section, we are interested in finding a minimal penalty $V(h)$, beyond which the procedure fails. Indeed, if $a$ and then $V(h)$ is too small, 
the minimization of the criterion amounts to minimize the bias, and then to choose the smallest possible bandwidth. This leads to the worst estimator and 
the risk explodes. 

In the following result $h_{\min}$ denotes the smallest bandwidth in $\H$ and is of order $1/n$.

\begin{thm} \label{penaliteminimale}
Assume  that $f$ is bounded. Choose $\H=\{e^{-k}, \lceil 2\log\log n\rceil \leq k \leq \lfloor\log n\rfloor\}$ as a  set of bandwidths.
Consider for $K$ the Gaussian kernel, the rectangular kernel, the Epanechnikov kernel or the biweight kernel.
If $a<1$, then there exists $C>0$ (depending on $f$, $a$, $K$) such that, 
for $n$ large enough (depending on $f$ and $K$), the selected bandwidth $\hat h$ defined by \eqref{defBV} and \eqref{defhath} satisfies
$$\P(\hat h \geq 3h_{\min})\leq C (\log n)^2\exp(-(\log n)^2/C)$$
i.e. $\hat h <3 h_{\min}$ with high probability.
Moreover $$\liminf_{n\to\infty}\E\|f-\hat f_{\hat h}\|^2>0.$$ 
\end{thm}

This theorem is proved in Section~\ref{preuvepenmin} for more general kernels and bandwith sets. Here we have simplified the conditions on $\mathcal{H}$ for the sake of readability. 
Actually the real condition on 
$\mathcal{H}$ for Theorem 3 is  that $E_\H=\min \{h/h' ; h\in \H, h'\in \H, h>h'\}$ does not depend on $n$ and is larger than 1. 
 It can be verified for the highlighted set $\H=\{e^{-k}, a_n \leq k \leq b_n\}$,
 but for $\H=\{c_n+d_nk, a_n \leq k \leq b_n\}$ as well. 
 
Mathematically, the proof of this result  relies on two main arguments. The first argument is probabilistic: roughly speaking concentration inequalities which allow to deal with expectations of the pairwise square distances between estimators instead of the square distances themselves. The other argument is analytical: it essentially relies on proper substitutes to Pythagoras' formula for kernel smoothing. 
The phase transition phenomenon is actually easier to highlight in a context for which we have the actual Pythagoras' identity at our disposal, see the discussion on projection estimators for Gaussian white noise model in \cite{nousArxiv}.
 
Theorem~\ref{penaliteminimale} ensures that the critical value for the parameter $a$ is 1. Beyond this value, the selected bandwidth $\hat h$ is of order $1/n$, which is very small 
(remember that for minimax study of a density with regularity $\alpha$, the optimal bandwidth is $n^{-1/(2\alpha+1)}$), then the risk cannot tend to 0. 

\bigskip

\section{Simulations}\label{simus}
  
  In this Section, we illustrate the role of tuning parameter $a$, the constant in the penalty term $V$. The aim is to observe the evolution of the risk 
  for various values of $a$. Is the critical value $a=1$ observable in practice? To do this, we simulate data $X_1, \dots, X_n$ for several densites $f$. 
  Next, for a grid of values for $a$, we compute the selected bandwidth $\hat h$, 
  the estimator $\hat f_{\hat h}$ and the integrated loss $\|\hat f_{\hat h}-f\|^2$.
  
\begin{figure}[h]
\includegraphics[scale=0.45]{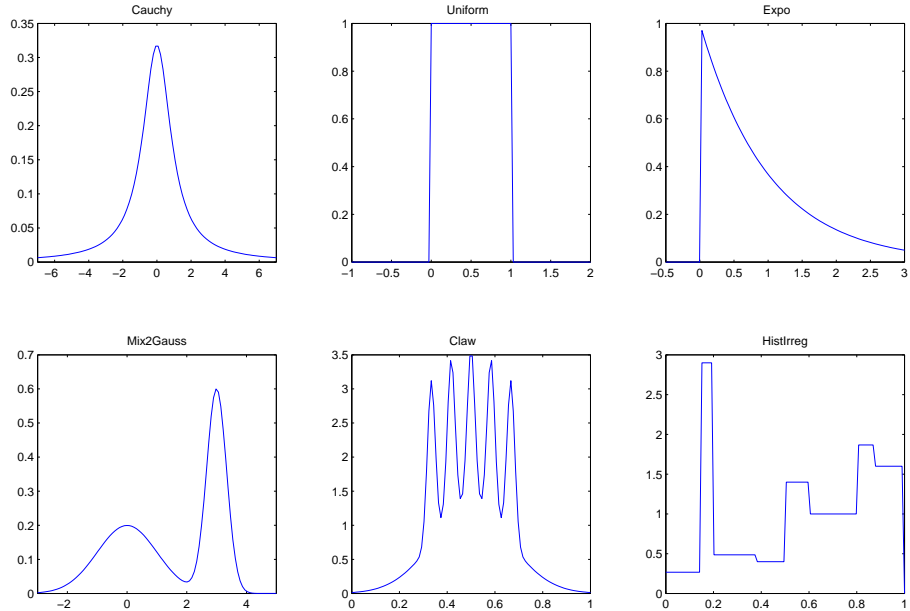}
\caption{Plots of true density $f$ for Examples 1--6}
\label{figdensites}
\end{figure}
  
  We consider the following examples, see Figure~\ref{figdensites}: 
  \begin{description}
  \item[Example 1]
   $f$ is the Cauchy density
    \item[Example 2]
   $f$ is the uniform density $\mathcal{U}(0,1)$
    \item[Example 3]
   $f$ is the exponential  density $\mathcal{E}(1)$ 
    \item[Example 4]
   $f$ is a mixture of two normal densities $\frac12\mathcal{N}(0,1)+\frac12\mathcal{N}(3,9)$
   \item[Example 5]
   $f$  is a mixture of  normal densities sometimes called Claw
   \item[Example 6]
   $f$ is a mixture of eight uniform densities
  \end{description}
  We implement the method for various kernels, but we only present results for Gaussian kernel, since the choice of kernel
  does not modify the results. 
  On the other hand, the method is sensitive to the choice of bandwidths set $\H$: here we use 
  $$\H=\{e^{-k}, 3 \leq k \leq 10\}\cup \{0.002+k\times0.02, 0\leq k \leq 24\}.$$
  Note that the theoretical conditions on the bandwidths are asymptotic. Then, they have no real sense in our simulations with given $n$. 
 In practice, this set must be rich enough for catching optimal bandwidths for a large class of densities, but small enough for the computation time.
 For our study, we choose equally distributed bandwidths for a good observation of the choice of $\hat h$, and we also add the set $\{e^{-k}, 3 \leq k \leq 10\}$ to have very small bandwidths avalaible, which are useful for irregular densities.

\begin{figure}[h]
\hspace{-2cm}
\begin{tabular}{ll}
\includegraphics[scale=0.45]{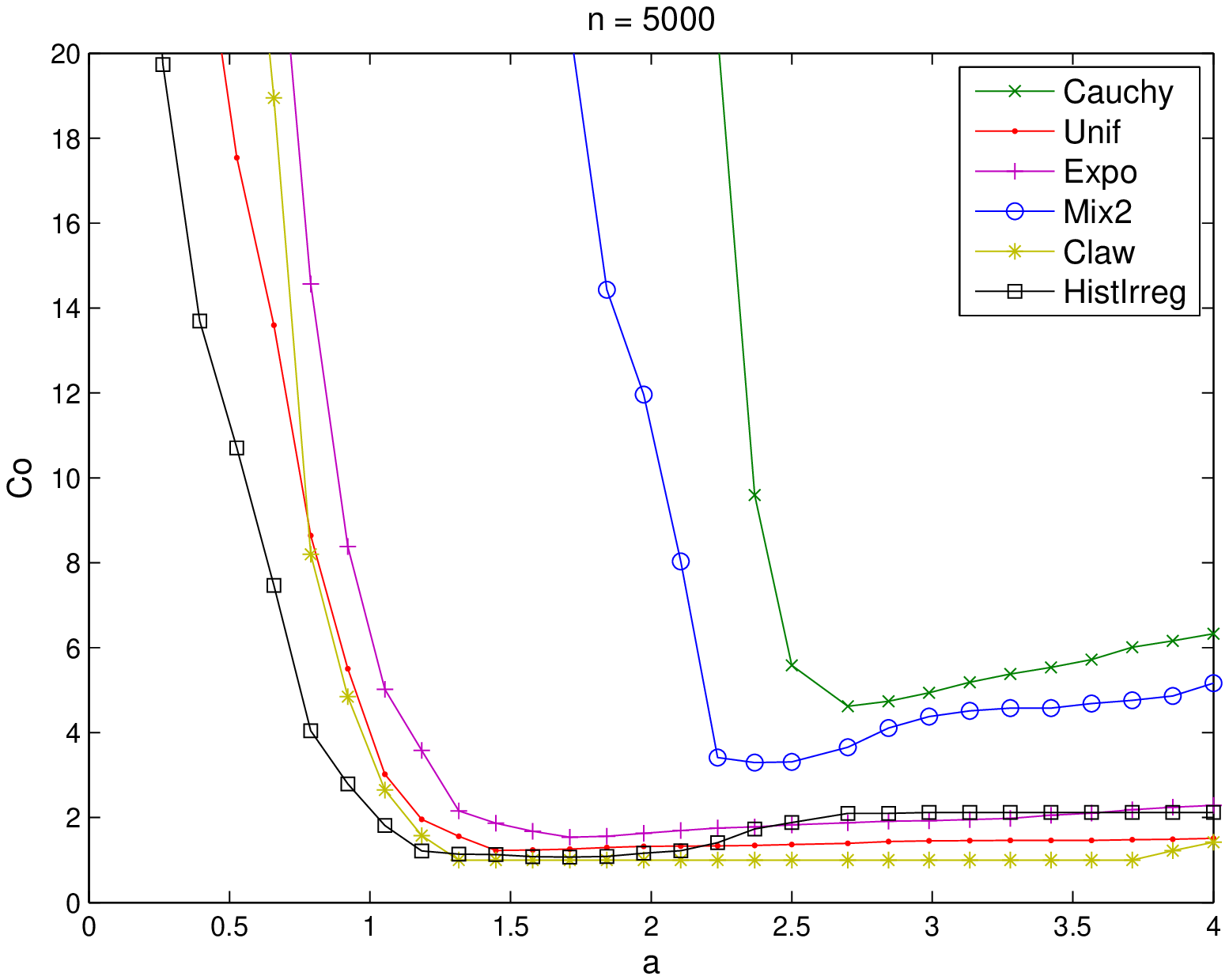}
&\includegraphics[scale=0.45]{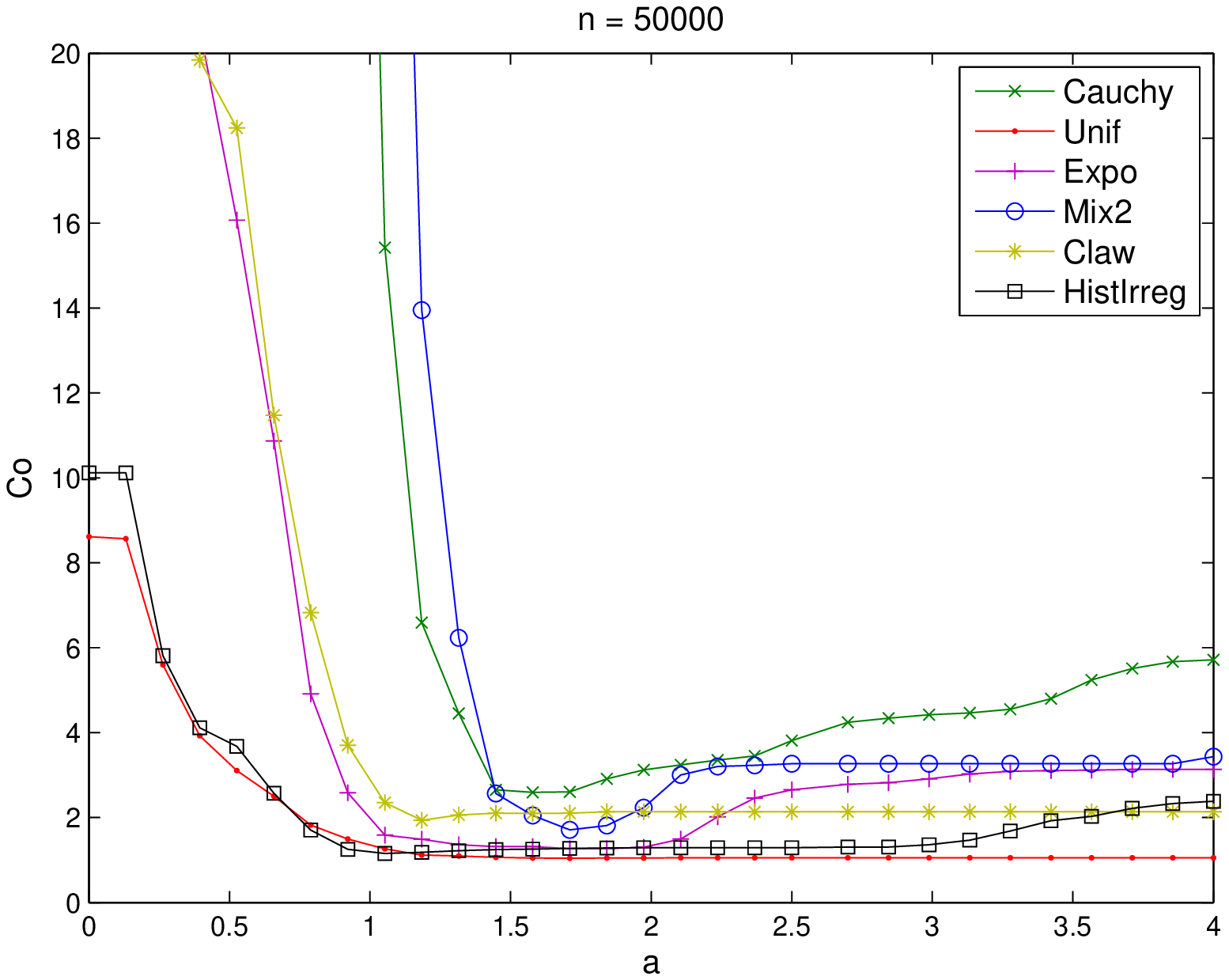}
\end{tabular}
\caption{Oracle constant $C_0$ as a function of $a$, for Examples 1--6}
\label{figCo}
\end{figure}

\begin{figure}[h]
\hspace{-2cm}
\begin{tabular}{ll}
\includegraphics[scale=0.45]{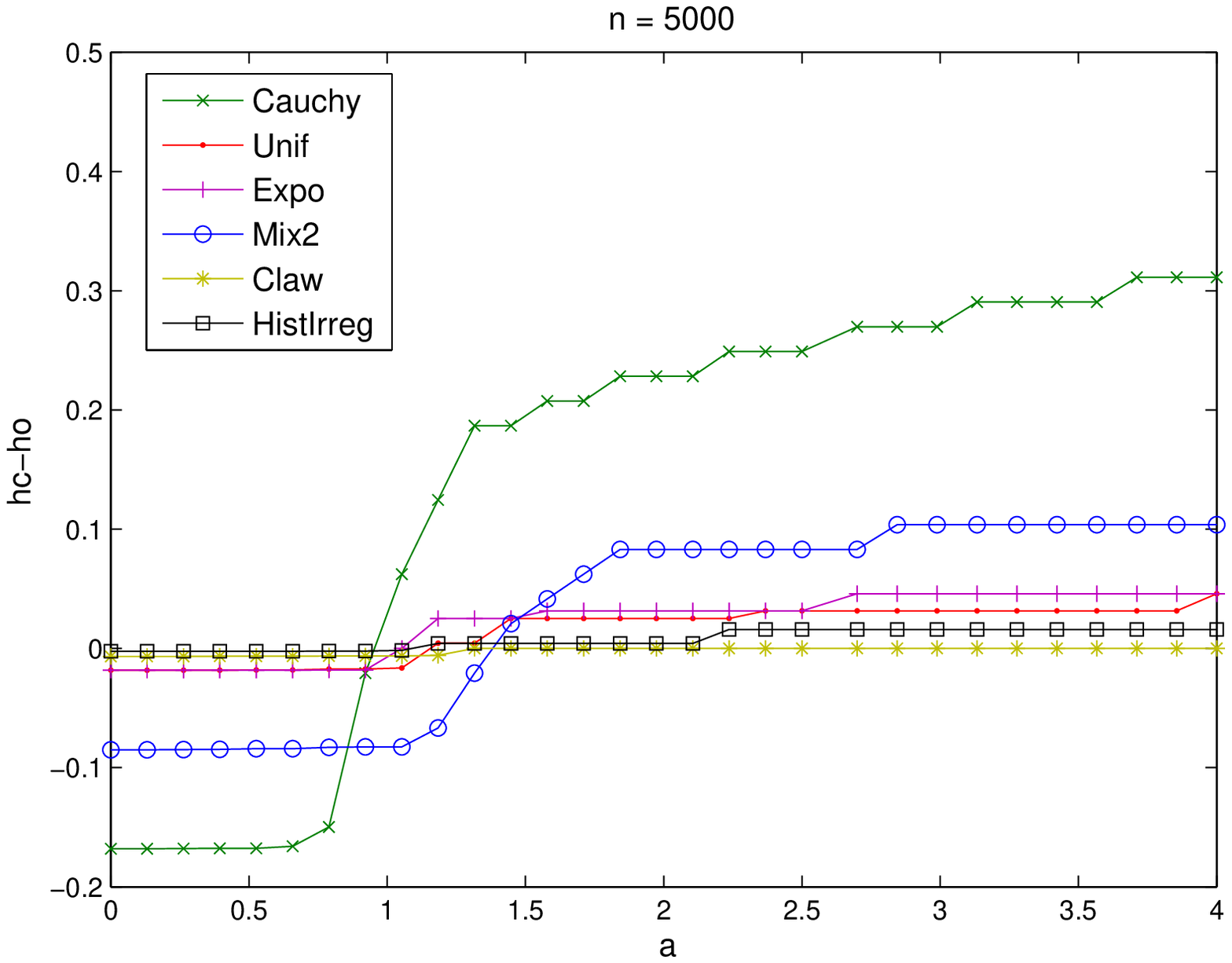}
&\includegraphics[scale=0.45]{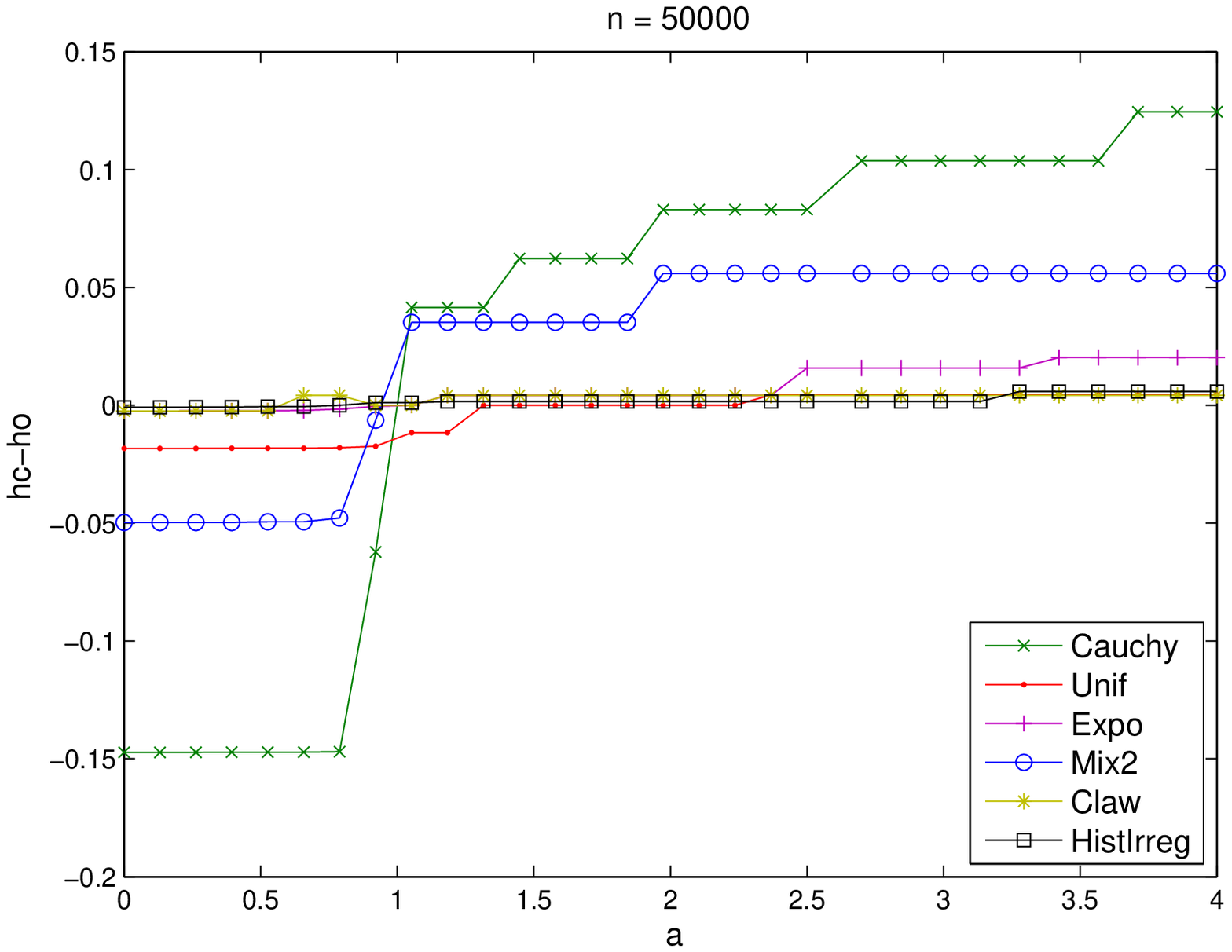}
\end{tabular}
\caption{$\hat h - h_0$ as a function of $a$, for Examples 1--6}
\label{fighc}
\end{figure}

 For $n=5000$ and $n=50000$, and  several values of $a$, the Figure~\ref{figCo} plots
 $$C_0=\tilde\E\frac{\|\hat f_{\hat h}-f\|^2}{\min_{h\in \H} \|\hat f_{h}-f\|^2 }$$
 where $\tilde\E$ means the empirical mean on $N=50$ experiments. 
 Thus smaller $C_0$ better the estimation. 
 Moreover, we also plot on Figure~\ref{fighc} the selected bandwidth compared to the optimal bandwidth in the selection (for $N=1$ experiment),
i.e. 
$$\hat h - h_0 \qquad\text{ where } \quad\|\hat f_{h_0}-f\|^2=\min_{h\in \H} \|\hat f_{h}-f\|^2.$$

We can observe that the risk (and then the oracle constant $C_0$) is very high for small values of $a$, as expected. Then it jumps to a small value, that indicates the method begins to work well. For too large values of $a$ the risk finally goes back up. 
Thus we observe in practice the transition phenomenon that was announced by the theory.
However, contrary to the theoretical results, the critical value may be not exactly at $a=1$, especially for small values of $n$.
As already mentioned above this is related to the asymptotic nature of the theoretical results that we have obtained.
For irregular densities (examples 2, 5, 6), the optimal bandwidth is very low, then it is consistent to observe a smaller jump for the bandwidth choice. 
However the jump does exist and this is the interesting point.
We can also observe that the optimal value for $a$ seems to be very close to the jump point. That may pose a problem of calibration and this is what we would like to discuss now.
  
\section{Discussion}

To calibrate the penalty $V$, we face two practical problems: 
first, the optimal value for $a$ seems to be extremely close to the minimal value; secondly, this latter value is not necessarily equal to the (asymptotic) theoretical value $a=1$.
In order to clearly separate the optimal value from the minimal, we propose to use some slightly different procedure, which depends on two possibly different penalty parameters instead of one as in the previous one.
\begin{eqnarray*}
&& B(h)=\sup_{h'\leq h}\left[\|\hat f_{h'}-\hat f_{h}\|^2-a\frac{\|K_{h'}\|^2}n\right]_+,\\
&& \hat{h}=\arg\min_{h\in\mathcal{H}}\left\{ B(h)+b\frac{\|K_{h}\|^2}n\right\} 
\end{eqnarray*}
with $b \neq a$.
Of course this procedure is merely the one that we have previously studied when $a=b$. 
Our belief is that taking $a$ and $b$ to be different leads to a better and more stable calibration. 
A good track for practical purpose seems to use the procedure of Section~\ref{sec:sup} to find $\hat a$ where there is a jump in the risk (in practice this jump can be detected on the selected bandwidths) and then to choose $b=2\hat a$. Once again, what is important for practical calibration of the penalty is not that the jump appears at  $a=1$ (this value should be considered as some "asymptopia" which is never achieved) but that the jump does exist so that it becomes possible to use the calibration strategy that we just described. Proving theoretical results for this procedure is another interesting issue related to optimality considerations for the penalty that we do not intend to address here.
\section{Proofs}\label{preuves}

\subsection{Proof of Proposition~\ref{oi}}\label{preuveoi}


The first step is to write, for some fixed $h\in \H$, 
$$\|\hat f_{\hat h}-f\|\leq \|\hat f_{\hat h}-\hat f_{h}\|+\|\hat f_{h}-f\|.$$
The last term can be splitted in $\|\hat f_{h}-f_h\|+\|f_{h}-f\|\leq \|\hat f_{h}-f_h\|+D(h)$.
Notice that for all $h'\leq h$, using \eqref{defBV}, 
$\|\hat f_{h'}-\hat f_{h}\|^2\leq B(h)+V(h')$, which can be written , for all $h, h'$; 
$$\|\hat f_{h'}-\hat f_{h}\|^2\leq B(h\vee h')+V(h\wedge h')$$
where $h\vee h'=\max(h,h')$ and $h\wedge h'=\min (h,h')$.
Then, using \eqref{defhath}, 
$$\|\hat f_{\hat h}-\hat f_{h}\|^2\leq B(h\vee\hat h)+V(h\wedge  \hat h)\leq B(h)+V(h)+ \max (B(h), V(h)).$$
We obtain, for any $h\in \H$, 
$$\|\hat f_{\hat h}-f\|\leq \sqrt{ 2B(h)+2V(h)}+D(h)+\|\hat f_{h}-f_h\|.$$
Thus the heart of the proof is to control $B(h)=\sup_{h'\leq h}\left[\|\hat f_{h'}-\hat f_{h}\|^2-V(h')\right]_+$ by a bias term. 
First we center the variables and write
$$\|\hat f_{h'}-\hat f_{h}\|^2\leq (1+\varepsilon)\|\hat f_{h'}-f_{h'}-\hat f_{h}+f_h\|^2+(1+\varepsilon^{-1})\|f_{h'}-f_h\|^2,$$
with $\varepsilon$ some positive real to specified later.
Moreover $\|\hat f_{h'}-f_{h'}-\hat f_{h}+f_{h}\|=\sup_{t\in B}\nu(t)$ where $B$ is the unit ball in $L^2$
 and $$\nu(t)=\langle t, \hat f_{h'}-f_{h'}-\hat f_{h}+f_{h}\rangle
=\frac1n\sum_{i=1}^n g_t(X_i)-\E(g_t(X_i))$$
with
$$g_t(X)=\int (K_{h'}-K_h)(x-X)t(x)dx.$$

We shall now use the  concentration inequality stated in Lemma~\ref{conc}, with $\mathcal F$ a countable set in $B$
such that  $\sup_{t\in {\mathcal F}}\nu(t)=\sup_{t\in { B}}\nu(t)$
(this equality is true for any dense subset of $B$ for the $L^2$ topology, since $\nu$ is continuous).
To apply result \eqref{concsup}, we need to compute $b$, $H$ and $v$. 
\begin{itemize}
\item For all $y\in\R$, since $t\in B$,
$$|g_t(y)|=|\int (K_{h'}-K_h)(x-y)t(x)dx|\leq \|K_{h'}-K_h\| \|t\|\leq \|K_{h'}-K_h\|\leq \|K_{h'}\|$$
so that $b=\|K_{h'}\| .$ We used assumption {\bf (K0)} which implies, for $h'\leq h$, $\|K_{h'}-K_h\|^2\leq \|K_{h'}\|^2-\|K_h\|^2\leq  \|K_{h'}\|^2$.
\item Jensen's inequality gives $H^2\leq \E(\sup_{t\in {\mathcal F}}\nu^2(t))$.
Now 
\begin{eqnarray}
\sup_{t\in \mathcal F}\nu^2(t)&=&\|\hat f_{h'}-f_{h'}-\hat f_{h}+f_{h}\|^2\nonumber\\
&=&\|\frac 1n \sum_{i=1}^n (K_{h'}-K_h)(x-X_i)-\E((K_{h'}-K_h)(x-X_i))\|^2\nonumber\\
\E(\sup_{t\in \mathcal F}\nu^2(t))&=&\int \var (\frac1n \sum_{i=1}^n (K_{h'}-K_h)(x-X_i))dx\nonumber\\
&=&\frac1n \int \var ( (K_{h'}-K_h)(x-X_1))dx\label{varnu}\\
&\leq & \frac1n \int \E ( (K_{h'}-K_h)^2(x-X_1))dx\\
&\leq & \frac1n \|K_{h'}-K_h\|^2 \leq  \frac1n \|K_{h'}\|^2\nonumber
\end{eqnarray}
Then $H^2\leq n^{-1}\|K_{h'}\|^2$.

\item For the variance term, let us write
 \begin{eqnarray*}
 \var (g_t(X_1)) &\leq & \E\left[\left( \int (K_{h'}-K_h)(x-X)t(x)dx)\right)^2\right]\\
& \leq&\E\left[\int |K_{h'}-K_h|(x-X)dx\right]\E\left[ \int |K_{h'}-K_h|(x-X)t^2(x)dx\right]\\
 &\leq &\|K_{h'}-K_h\|_1^2\|f\|_\infty\|t\|^2\leq 4\|K\|_1^2\|f\|_\infty\|t\|^2
\end{eqnarray*}
since $\|K_{h'}-K_h\|_1\leq 2\|K\|_1$. Then $v= 4\|K\|_1^2\|f\|_\infty=4\|f\|_\infty.$
\end{itemize}

Finally, using \eqref{concsup}, with probability larger than $1-\sum_{h'< h}\max(e^{-\frac{\varepsilon^2\wedge \varepsilon}{24}\sqrt{n}},
e^{-\frac{\varepsilon^2\|K\|^2}{24\|f\|_\infty}\frac{1}{h'}})$
   \begin{eqnarray*}\label{conc2}
  \forall  h'\leq h \in \mathcal{H}\qquad \|\hat f_{h'}-f_{h'}-\hat f_{h}+f_h\|\leq 
(1+\varepsilon)\frac{\|K_{h'}\|}{\sqrt{n}}
  \end{eqnarray*}
  where we choose $\varepsilon$ such that $a\geq (1+\varepsilon)^3$.
Then, with probability larger than $1-\sum_{h\in \H}\sum_{h'\leq h}\max(e^{-\frac{\varepsilon^2\wedge \varepsilon}{24}\sqrt{n}},
e^{-\frac{\varepsilon^2\|K\|^2}{24\|f\|_\infty}\frac{1}{h'}})$
for any $h$, 
 $$B(h)\leq (1+\varepsilon^{-1}) D(h)^2$$
 In the same way, choosing  $0<\epsilon \leq \sqrt{a}-1$, we can prove that, 
 with probability $1-\sum_{h\in \H}\max(e^{-\frac{\epsilon^2\wedge \epsilon}{24}\sqrt{n}},
e^{-\frac{\epsilon^2\|K\|^2}{6\|f\|_\infty}\frac{1}{h}})$,
for any $h$, 
 $$\|\hat f_{h}-f_h\|\leq (1+\epsilon) \frac{\|K_{h}\|}{\sqrt{n}}\leq \sqrt{V(h)}$$
Finally, with high probability,  
\begin{eqnarray*}
\|\hat f_{\hat h}-f\|&\leq &\sqrt{2(1+\varepsilon^{-1}) D(h)^2+2V(h)}+D(h)+\sqrt{V(h)}\\
&\leq &(\sqrt{2(1+\varepsilon^{-1})}+1)\left(D(h)+\sqrt{V(h)}\right)
\end{eqnarray*}
To conclude we choose $\varepsilon=\epsilon=a^{1/3}-1$.
Regarding the second result, note that 
the rough bound $\|\hat f_h\|^2\leq \|K_h\|^2\leq \|K\|^2/h_{\min}$ is valid for all $h$. Then,
denoting $A$ the set on which the previous oracle inequality is verified,  
$$\E\|\hat f_{\hat h}-f\|^2\leq \E\|\hat f_{\hat h}-f\|^2\1_A+2(\|f\|^2+\|K\|^2/h_{\min})\P(A^c)$$
with $$\P(A^c)\leq 2\sum_{h,h'\in \H}\max(e^{-c\sqrt{n}},e^{-\frac{c}{h'}})
\leq 2|\H|^2e^{-c/{h_{\max}}}$$

\hfill$\blacksquare$\\

\subsection{Proof of  Theorem~\ref{penaliteminimale}}\label{preuvepenmin}

We shall prove a more general version of the theorem, where several bandwidths sets $\H$ and kernels $K$ are possible. 
We denote $\Crit(h):=B(h)+V(h)$ and 
$E_\H=\min \{h/h' ; h\in \H, h'\in \H, h>h'\}$.
We assume that $E_\H$ does not depend on $n$ and is larger than 1 ($\H=\{e^{-k}, a_n \leq k \leq b_n\}$ suits with $E_\H=e$).
 Let us define $$\phi(x)=\|K\|^{-2}\|K-K_x\|^2=1+\frac{1}x-2 \frac{\langle K, K(x.) \rangle}{\|K\|^2}.$$ 
We assume that the kernel $K$ satisfies : 
\begin{description}
\item[(K1)] the function $\phi$ is bounded from below over $[E_\H,+\infty)$,
\item[(K2)] for $0<\mu<1$, the function $\phi(x)-\frac{\mu}x$ tends to $+\infty$ when $x\to 0$ and is decreasing in some neighborhood of $0$,
\item[(K3)] for $0<\mu<1$, the function  $\phi(x)+\frac{\mu}x$ is increasing for  $x\geq 2$.
\end{description}

 These assumptions are mild, as shown in the following Lemma, proved in Section~\ref{preuveL4}. 
 
 \begin{lemma}\label{noyaux}
 The following kernels satisfy assumptions {\bf (K0--K3)}:
 \begin{itemize}
  \item[a - ] Gaussian kernel: $K(x)=e^{-x^2/2}/\sqrt{2\pi}$
  \item[b - ] Rectangular kernel: $K(x)=\1_{[-1,1]}(x)/2$
  \item[c - ] Epanechnikov kernel: $K(x)=(3/4)(1-x^2)\1_{[-1,1]}(x)$
  \item[d - ] Biweight kernel: $K(x)=(15/16)(1-x^2)^2\1_{[-1,1]}(x)$
  \end{itemize}
 \end{lemma}

 The general result is:

\begin{thm} Assume {\bf (K0--K3)} and that $f$ is bounded. 
Assume that $E_\H$ does not depend on $n$ and $h_{\max}\to 0$ when $n\to\infty$. 
We also assume that there exist $\theta_1<\theta_2$ reals such that $\theta_2\geq 2$, $\theta_1.h_{\min}\in \H$ and 
$\phi(\theta_2)-\phi(\theta_1)\geq 1/\theta_1-1/\theta_2$. 

Then,  if $a<1$, there exists $C=C(\|f\|_\infty)>0$ such that, for $n$ large enough (depending on $f, \H, K$), 
$$ \P(\hat h \geq \theta_2h_{\min})\leq \sum_{h\in \H}\sum_{h'< h}\max(e^{-C{\varepsilon^2}{\sqrt n}},e^{-C\varepsilon^2\|K_{h'}-K_h\|^2})$$
where $\varepsilon<1-a^{1/3}$. If $\H=\{e^{-k}, a_n \leq k \leq b_n\}$ and the kernel is Gaussian, rectangular, Epanechnikov or biweight, $\theta_1=e$ and $\theta_2=3$ work.
\end{thm}

This results implies  Theorem~\ref{penaliteminimale}, since under {\bf (K1)}, $\|K_{h'}-K_h\|^2
=\frac{\|K\|^2}{h'}\phi(h/h')
\geq (\min_{E_\H}\phi) \frac{\|K\|^2}{h'} $ as soon as $h>  h'$,
so that 
$$\sum_{h\in \H}\sum_{h'< h}e^{-C\|K_{h'}-K_h\|^2}\leq |\H|^2e^{-C'/h_{\max}}.$$

Let $\varepsilon\in (0,1)$ such that $a<(1-\varepsilon)^3$ and 
\begin{equation}\label{theta}
\varepsilon^3+3\varepsilon<\frac{\phi(\theta_2)-\phi(\theta_1)-a/\theta_1+a/\theta_2}{\phi(\theta_2)+\phi(\theta_1)}
\end{equation}(possible since $a<1\leq (\phi(\theta_2)-\phi(\theta_1))/(1/\theta_1-1/\theta_2)$).
Let us decompose
$$\hat f_{h'}-\hat f_{h}=(\hat f_{h'}-f_{h'}-\hat f_{h}+f_{h})+( f_{h'}- f_{h})=S(h,h')+( f_{h'}- f_{h})$$
 with $$S(h,h')=\frac1n\sum_{i=1}^n (K_{h'}-K_h)(x-X_i)-\E((K_{h'}-K_h)(x-X_i))$$
and the bias term
$\|f_{h'}- f_{h}\|\leq \sup_{h'\leq h}\| K_{h'}*f- K_{h}*f\|=D(h)$.
First write
$$(1-\varepsilon) \|S(h,h')\|^2-\left(\frac1\varepsilon-1\right)D(h)^2\leq \|\hat f_{h'}-\hat f_{h}\|^2\leq (1+\varepsilon) \|S(h,h')\|^2+\left(1+\frac1\varepsilon\right)D(h)^2$$
Now we shall prove that with high probability 
$$(1-\varepsilon)\frac{\|K_{h'}-K_h\|}{\sqrt{n}}\leq\|S(h,h')\|\leq (1+\varepsilon)\frac{\|K_{h'}-K_h\|}{\sqrt{n}}.$$

First, we can prove as in Section~\ref{sec:sup} that for all $h'<h$
\begin{eqnarray*}
&&\P\left(\|S(h,h')\|\geq (1+\varepsilon)\frac{\|K_{h'}-K_h\|}{\sqrt{n}}\right)\\
&&\leq \max\left(\exp\left(-\frac{\varepsilon^2\wedge \varepsilon}{24}\sqrt{n}\right),
\exp\left(-\frac{\varepsilon^2}{24\|f\|_\infty}\|K_{h'}-K_h\|^2\right)\right).
\end{eqnarray*}
Next, we shall use \eqref{concinf} in Lemma~\ref{conc} in order to lowerbound $\|S(h,h')\|$.
Recall that $\|S(h,h')\|=\sup_{t\in B}\nu(t)$ where $B$ is the unit ball in $L^2$
 and $\nu(t)=\frac1n\sum_{i=1}^n g_t(X_i)-\E(g_t(X_i))$
with $g_t(X)=\int (K_{h'}-K_h)(x-X)t(x)dx.$ 
With notations of Lemma~\ref{conc}, we have $b=\|K_{h'}-K_h\|$, $H'^2=n^{-1}\|K_{h'}-K_h\|^2$ and $v=4\|K\|_1^2\|f\|_\infty$. It remains to lowerbound $H$.
First, remark, that \eqref{varnu} provides $n\E(\sup_{t\in B}\nu^2(t))=\|K_{h'}-K_h\|^2- \| (K_{h'}-K_h)*f\|^2$. Next,
using \eqref{varsup}
$$\E (\sup_{t\in B}\nu^2(t))\leq \frac vn+4\frac{bH}n+H^2\leq \frac vn+\left(H+\frac{2b}n\right)^2.$$
Then 
$$ n\left(H+\frac{2b}n\right)^2\geq n\E (\sup_{t\in B}\nu^2(t))- v=\|K_{h'}-K_h\|^2- \| (K_{h'}-K_h)*f\|^2-4\|K\|_1^2\|f\|_\infty$$
which implies
$$ \sqrt{n}\left(H+\frac{2b}n\right)\geq \sqrt{ \|K_{h'}-K_h\|^2 -4\|K\|_1^2(\|f\|_\infty+\|f\|^2)}.$$
Since $b=\|K_{h'}-K_h\|$,
$$H\geq \sqrt{\frac{\|K_{h'}-K_h\|^2 -4\|K\|_1^2(\|f\|_\infty+\|f\|^2)}{n}}-\frac{2\|K_{h'}-K_h\| }n$$
Now, for $h'< h$
$$H\geq \frac{\|K_{h'}-K_h\| }{\sqrt{n}}\left(\sqrt{1-\frac{ 4\|K\|_1^2(\|f\|_\infty+\|f\|^2)}{\|K_{h'}-K_h\|^2}}-\frac2{\sqrt{n}}\right)$$
so
$$H-\frac{\varepsilon}{3} H'\geq H' \left(\sqrt{1-\frac{ 4\|K\|_1^2(\|f\|_\infty+\|f\|^2)}{\|K_{h'}-K_h\|^2}}-\frac2{\sqrt{n}}-\frac{\varepsilon}{3}\right).$$

From {\bf (K1)}, $\|K_{h'}-K_h\|^2=\frac{\|K\|^2}{h'}\phi(h/h')
\geq (\min_{E_\H}\phi) \frac{\|K\|^2}{h'} \geq\frac{C}{h_{\max}}\to \infty$ and, in consequence, for $n$ large enough 
$$H-\frac{\varepsilon}{3} H'\geq H' \left(1-\varepsilon\right).$$
Thus for $n$ large enough
\begin{eqnarray}\label{aa}
&&\P\left(\|S(h,h')\|\leq (1-\varepsilon)\frac{\|K_{h'}-K_h\|}{\sqrt{n}}\right)\\
&&\leq\max\left(\exp\left(-\frac{\varepsilon^2\wedge (3\varepsilon)}{24\times 9}\sqrt{n}\right),
\exp\left(-\frac{\varepsilon^2}{24\times 9\|f\|_\infty}\|K_{h'}-K_h\|^2\right)\right)\nonumber
\end{eqnarray}
Let $\delta(h,h)=0$ and, if $h\neq h'$,
$$\delta(h,h')=2\max\left(\exp\left(-\frac{\varepsilon^2\wedge (3\varepsilon)}{24\times 9}\sqrt{n}\right),
\exp\left(-\frac{\varepsilon^2}{24\times 9\|f\|_\infty}\|K_{h'}-K_h\|^2\right)\right).$$
We just proved that for $n$ large enough, with probability larger than $1-\delta(h,h')$
$$(1-\varepsilon)^2\frac{\|K_{h'}-K_h\|^2}{{n}}\leq\|S(h,h')\|^2\leq (1+\varepsilon)^2\frac{\|K_{h'}-K_h\|^2}{{n}}.$$

Next, with probability larger than $1-\sum_{h'\leq h}\delta(h,h')$
$$\begin{cases}
 B(h)\geq \sup_{h'\leq h}\left[(1-\varepsilon)^3\frac{\|K_{h'}-K_h\|^2}{{n}}-a\frac{\|K_{h'}\|^2}n\right]_+-\left(\frac1\varepsilon-1\right)D(h)^2\\
 B(h)  \leq
  \sup_{h'\leq h}\left[(1+\varepsilon)^3\frac{\|K_{h'}-K_h\|^2}{{n}}-a\frac{\|K_{h'}\|^2}n\right]_++\left(1+\frac1\varepsilon\right)D(h)^2
  \end{cases}
$$
But, if $h_{\min}$ small enough, for $\lambda>a$ 
$$ \sup_{h'\leq h}\left[\lambda\frac{\|K_{h'}-K_h\|^2}{{n}}-a\frac{\|K_{h'}\|^2}n\right]_+
 =\lambda\frac{\|K_{h_{\min}}-K_h\|^2}{{n}}-a\frac{\|K_{h_{\min}}\|^2}n$$
Indeed, for $x=h'/h\leq 1$
\begin{eqnarray*}
\lambda\frac{\|K_{h'}-K_h\|^2}{{n}}-a\frac{\|K_{h'}\|^2}n&=&\lambda\frac{\|K\|^2}{nh}\left(1+\frac{1-a/\lambda}x-2 \frac{\langle K, K(x.) \rangle}{\|K\|^2}\right)\\
&=&\lambda\frac{\|K\|^2}{nh}\left(\phi(x)-\frac{a/\lambda}x\right)
\end{eqnarray*}
 and the function $\phi(x)-\frac{a/\lambda}x$ tends to $+\infty$ when $x\to 0$ and is decreasing in some neighborhood of $0$ (assumption {\bf (K2)}).
 %
Then with probability larger than $1-\sum_h\sum_{h'\leq h}\delta(h,h')$, for all $h$
$$\begin{cases}
\Crit(h)\geq \frac{\|K\|^2}{nh_{\min}}\left[-a+(1-\varepsilon)^3\phi(h/h_{\min})+\frac a{h/h_{\min}}\right]-\left(\frac1\varepsilon-1\right)D(h)^2\\
\Crit(h)\leq \frac{\|K\|^2}{nh_{\min}}\left[-a+(1+\varepsilon)^3\phi(h/h_{\min})+\frac a{h/h_{\min}}\right]+\left(1+\frac1\varepsilon\right)D(h)^2
\end{cases}$$

In particular, for $h=\theta_1 h_{\min}$, 
\begin{equation}\label{theta1}
\Crit(\theta_1 h_{\min})\leq \frac{\|K\|^2}{nh_{\min}}\left[-a+(1+\varepsilon)^3\phi(\theta_1)+\frac a{\theta_1}\right]+\left(1+\frac1\varepsilon\right)\sup_hD(h)^2.
\end{equation}

Moreover, since $a<(1-\varepsilon)^3$, $(1-\varepsilon)^3\phi(x)+\frac {a}{x}$ is increasing for  $x\geq 2$ (assumption {\bf (K3)}). 
This implies that 
\begin{equation}\label{theta2}
\forall h\geq \theta_2 h_{\min},\quad
\Crit(h)\geq \frac{\|K\|^2}{nh_{\min}}\left[-a+(1-\varepsilon)^3\phi(\theta_2)+\frac a{\theta_2}\right]-\left(\frac1\varepsilon-1\right)\sup_hD(h)^2.
\end{equation}%

Since $(K_h)$ is an approximation to the identity,  $\|f-K_h*f\|$ tends to $0$ when $h$ tends to 0. This implies that 
$D(h)\leq 2\sup_{h'\leq h}\| f- K_{h'}*f\|$ tends to $0$ and $\sup_{h \in \H}D(h)$ tends to $0$, as soon as $h_{\max}$ tends to 0.
Now \eqref{theta} leads to $\Delta:=(1-\varepsilon)^3\phi(\theta_2)+\frac a{\theta_2}-(1+\varepsilon)^3\phi(\theta_1)-\frac a{\theta_1}>0$.
Then, 
for $n$ large enough, $(2/\varepsilon)\sup_hD(h)^2< \frac{\|K\|^2}{nh_{\min}} \Delta$ so that 
\begin{equation}\label{ss}
\begin{split}
&\frac{\|K\|^2}{nh_{\min}}\left[-a+(1+\varepsilon)^3\phi(\theta_1)+\frac a{\theta_1}\right]+\left(1+\frac1\varepsilon\right)\sup_hD(h)^2\\
&< \frac{\|K\|^2}{nh_{\min}}\left[-a+(1-\varepsilon)^3\phi(\theta_2)+\frac a{\theta_2}\right]-\left(\frac1\varepsilon-1\right)\sup_hD(h)^2
\end{split}
\end{equation}
Finally, combining \eqref{theta1} and \eqref{theta2} and \eqref{ss} gives $\hat h< \theta_2 h_{\min}$ with probability larger than $1-\sum_h\sum_{h'\leq h}\delta(h,h')$.

Let us now prove the second part of Theorem~\ref{penaliteminimale}, that is the lower bound on the risk.
Let $A_n=\{\hat h\leq 3h_{\min}\}$ and $B_n=\cap_{h\in \H}\{\|f_{ h}-\hat f_{h}\|\geq \frac12\frac{\|K_h\|}{\sqrt{n}}\}$.
We have just proved that $\P(A_n^c)\leq C (\log n)^2\exp(-(\log n)^2/C)$. In the same way that \eqref{aa}, we can write for $n$ large enough 
\begin{eqnarray*}
&&\P\left(\|f_{ h}-\hat f_{h}\|\leq (1-\varepsilon)\frac{\|K_h\|}{\sqrt{n}}\right)
\\ && \leq\max\left(\exp\left(-\frac{\varepsilon^2\wedge (3\varepsilon)}{24\times 9}\sqrt{n}\right),
\exp\left(-\frac{\varepsilon^2}{6\times 9\|f\|_\infty}\|K_h\|^2\right)\right)
\end{eqnarray*}
which implies $\P(B_n^c)\leq C' (\log n)\exp(-(\log n)^2/C')$ and then
$$\P(A_n\cap B_n)\geq 1 - o(1).$$
Then we can write 
\begin{eqnarray*}
\|f-\hat f_{\hat h}\|&\geq &\|f_{\hat h}-\hat f_{\hat h}\|\1_{A_n\cap B_n}-\|f-f_{\hat h}\|\\
& \geq & \min_{h\leq 3h_{\min}}\|f_{ h}-\hat f_{h}\|\1_{A_n\cap B_n}-\max_{h}\|f-f_{h}\|\\
& \geq & \min_{h\leq 3h_{\min}}\frac12\frac{\|K_h\|}{\sqrt{n}}\1_{A_n\cap B_n}-\max_{h}\|f-f_{h}\|\\
& \geq & \frac{\|K\|}{2\sqrt{3}}\frac{1}{\sqrt{nh_{\min}}}\1_{A_n\cap B_n}-\max_{h}\|f-f_{h}\|
\end{eqnarray*}
But $\max_{h}\|f-f_{h}\|\to 0$ (since $h_{\max} \to 0$), and $nh_{\min}\to 1$ when $n\to\infty$. Hence
\begin{eqnarray*}
\E\|f-\hat f_{\hat h}\|& \geq & \frac{\|K\|}{2\sqrt{3}}\frac{\P(A_n\cap B_n)}{\sqrt{1+o(1)}}-o(1) 
\end{eqnarray*}
which proves that  $\E\|f-\hat f_{\hat h}\| \geq \frac{\|K\|}{4\sqrt{3}}$ for $n$ large enough.

\hfill$\blacksquare$\\

\subsection{Proof of  Lemma~\ref{noyaux}}\label{preuveL4}

To prove Lemma~\ref{noyaux}, it is sufficient to do computations on integrals. We obtain:
\begin{itemize}
  \item[a - ] if $K$ is the Gaussian kernel,  $$\frac{\langle K, K(x.) \rangle}{\|K\|^2}=\sqrt{\frac{2}{1+x^2}}.$$
  
  \item[b - ] if $K$ is the rectangular kernel,  $$\frac{\langle K, K(x.) \rangle}{\|K\|^2}=\frac{1}{x}\wedge 1.$$
  
  \item[c - ] if $K$ is the Epanechnikov kernel,   $$\frac{\langle K, K(x.) \rangle}{\|K\|^2}=
 \frac{5}{4}  \left[\left(\frac{1}{x}\wedge 1\right)-\frac{x^2}{5}\left(\frac{1}{x}\wedge 1\right)^5\right].$$

 \item[d - ] if $K$ is the biweight kernel: $$\frac{\langle K, K(x.) \rangle}{\|K\|^2}=
 \frac{1}{16}  \left[21\left(\frac{1}{x}\wedge 1\right)-6{x^2}\left(\frac{1}{x}\wedge 1\right)^5+{x^4}\left(\frac{1}{x}\wedge 1\right)^9\right] .$$
 
These formulas permit to verify all the assumptions. 
  
  \end{itemize}

\hfill$\blacksquare$\\

\section*{References}

\bibliographystyle{elsarticle-harv} 
\bibliography{biblio}

\end{document}